\documentclass[letterpaper, 10 pt, conference]{ieeeconf}  

\IEEEoverridecommandlockouts                              

\usepackage{cite}
\usepackage{amsmath,amssymb,amsfonts}
\usepackage{graphicx}
\usepackage{textcomp}

\usepackage{graphics} 
\usepackage{epsfig} 
\usepackage{mathptmx} 
\usepackage{times} 
\usepackage{amsbsy}
\usepackage{dsfont}
\usepackage{verbatim}
\usepackage{algorithm}
\usepackage[colorlinks=true,urlcolor=blue]{hyperref}
\usepackage{algpseudocode}
\usepackage{url}
\usepackage{xcolor}
\graphicspath{{./images/}}

\DeclareMathOperator*{\argmin}{arg\,min}

\newcommand{\R}{\mathbb{R}}
\newcommand{\B}[1]{\mathbf{#1}}

\newcommand{\innerprod}[2]{\left\langle #1,#2 \right\rangle}
\newcommand{\abs}[1]{\left \vert #1 \right \vert}

\algrenewcommand\algorithmicindent{1em}

\begin{document}

\title{
A Hopf-Lax Type Formula for Multi-Agent Path Planning with Pattern Coordination 
}
\author{Christian Parkinson and Adan Baca
\thanks{*This work was supported by the NSF through the Research Training Group on Applied Mathematics and Statistics for Data Driven Discovery at the University of Arizona (DMS-1937229).}
\thanks{Christian Parkinson is an assistant professor in the Department of Mathematics and Department of Computational Mathematics, Science and Engineering, Michigan State University, East Lansing, MI, USA, 48824
        {\tt\small chparkin@msu.edu}}%
\thanks{Adan Baca is a student in the Department of Mathematics and Department of Computer Science, University of Arizona, Tucson, AZ, USA, 85721
        {\tt\small adanbaca@arizona.edu}}%
}

\maketitle

\begin{abstract}
We present an algorithm for a multi-agent path planning problem with pattern coordination based on dynamic programming and a Hamilton-Jacobi-Bellman equation. This falls broadly into the class of partial differential equation (PDE) based optimal path planning methods, which give a black-box-free alternative to machine learning hierarchies. Due to the high-dimensional state space of multi-agent planning problems, grid-based methods for PDE which suffer from the curse of dimensionality are infeasible, so we instead develop grid-free numerical methods based on variational Hopf-Lax type representations of solutions to Hamilton-Jacobi Equations. Our formulation is amenable to nonlinear dynamics and heterogeneous agents. We apply our method to synthetic examples wherein agents navigate around obstacles while attempting to maintain a prespecified formation, though with small changes it is likely applicable to much larger classes of problems. 
\end{abstract}


\thispagestyle{empty}

\section{INTRODUCTION}

The collaborative path planning problem involves multiple agents that share information navigating a domain while attempting to accomplish a set of tasks. This is a classical problem in automated path planning \cite{Coord1, Coord2} which is still of acute interest today \cite{CoordRev}. In particular, we are interested in the problem of multiple, possibly heterogeneous agents navigating while keeping some prespecified geometric formation. This problem has applications in terrain mapping \cite{CoordAirGround}, optimal transportation of multiple objects \cite{Coord4}, and wilderness search and rescue \cite{Coord5}.

Many modern planning algorithms rely on deep learning architectures \cite{DL1,DL2}, primarily because they are applicable to high-dimensional problems and able to handle complex constraints. However, learning algorithms like these can suffer from lack of robustness or explainability. While ameliorating these concerns is an active area of research \cite{Robust1,Robust2}, methods rooted in optimal control and partial differential equations (PDE) are often black-box-free and come with theoretical guarantees, making them a fully explainable alternative. 

PDE-based path planning methods have been used in simple models of autonomous vehicles \cite{ParkinsonCar1,ParkinsonCar2,TakeiTsai2}, human navigation \cite{Parkinson,Parkinson2}, environmental crime \cite{Arnold,Cartee2,Chen}, animal foraging in the presence of predation \cite{GeeFear}, and driving with traffic signal uncertainty \cite{Gaspard}, to list a few. These methods employ a dynamic programming \cite{Bellman} approach to optimal control wherein a value function encoding all necessary information regarding optimal travel is realized as the solution to a Hamilton-Jacobi-Bellman (HJB) equation. Classical numerical algorithms for approximating solutions to these equations can be broadly grouped into fast-marching schemes \cite{Tsitsiklis,SethVlad1} and fast-sweeping schemes \cite{FS1,FS4} each of which rely on discretizing a spatial domain and approximating derivatives via finite differences. Any such methods will suffer from the curse of dimensionality, wherein the number of grid nodes (and hence the computational complexity) scales exponentially with dimension, making them infeasible for high-dimensional systems or multi-agent problems. Accordingly, in recent years, there has been some effort devoted to developing grid-free numerical methods for HJB and related equations; see \cite{RecentAdvancesHJ} for a recent review of such methods. Many of these methods rely on Hopf-Lax type formulas, which provide the solution to an HJB equation at individual points in terms of an optimization problem. We recall the classical Hopf-Lax formula \cite[\S3.3]{Tran}, which says that under some mild conditions on the Hamiltonian $H$ and the initial function $g$, the solution of the Hamilton-Jacobi equation $u_t + H(\nabla u) = 0, u(x,0) = g(x)$ is \begin{equation} \label{eq:ClassicHL} u(x,t) = \inf_{x^*} \sup_{p^*} \Big\{g(x^*) + \langle x-x^*,p^*\rangle - tH(p^*) \Big\}. \end{equation} However, in path planning problems, any interesting model of motion will lead to a HJB equation with a Hamiltonian that is space and/or time dependent: $H = H(x,p,t)$ where $p$ is a proxy for $\nabla u$. The classical Hopf-Lax formula does not apply for such equations. In \cite{ParkinsonPolage,Lin,ParkinsonBoyle,Huynh}, a conjectural Hopf-Lax type formula is proposed and used to solve single agent path planning problems (as well as two-player differential games in \cite{Lin}).  Similar methods are used for a multi-agent problem in \cite{KirchnerCollabLinearized} through the dynamics are linearized which simplifies the computation.

Inspired by this, we present a control theoretic model for multi-agent control with geometric pattern formation, and a numerical method based on a discretized Hopf-Lax type formula which approximates the solution of the associated HJB equation.  We demonstrate the efficacy of our method with some sythetic examples of agents navigating a domain whle attemping to maintain a prespecified formation.  

\section{OPTIMAL CONTROL FORMULATION }

In this section, we present the optimal control problem that we will solve, and the corresponding Hamilton-Jacobi-Bellman equation. Throughout, we use boldface variables to denote state-space or control trajectories, and ordinary font to denote points in state space and individual control values. To shorten notation, we use brackets instead of tuples; for example $\{x_i\}$ is taken to mean $(x_1,x_2,\ldots,x_I)$.

We consider agents $\{\B x_i\}^I_{i=1}$ navigating according to dynamics \begin{equation}
\label{eq:dynamics}
\dot{\B x}_i = f_i(\B x_i,\B a_i,t), \,\,\,\,\, \B x_i(0) =\B x_{i,0}
\end{equation} where initial conditions $\B x_{i,0} \in \R^{n_i}$ and dynamics functions $f_{i}:\R^{n_i}\times \R^{m_i}\times [0,T] \to \R^{n_i}$ are given. Here $T > 0$ is the horizon time and $\B a_i : [0,T] \to \mathbb A_{i}\subset \R^{m_i}$ is the control for agent $i$. In particular, the agents could be heterogeneous so it is \emph{not} necessarily the case that the dimensions $n_i$ of the state spaces or $m_i$ of the control spaces are the same. We assume that the agents work collaboratively, or equivalently that there is a centralized controller choosing all controls. The goal of the controller is to minimize the cost functional which we decompose into two terms: \begin{equation} \label{eq:CostFunc}
\mathcal C[\{\B x_i\},\{\B a_i\}] = \int^T_0 w_1 \chi\Big(\{\B x_i(t)\}\Big) + w_2 \rho\Big(\{\B x_i(t)\}\Big)dt.
\end{equation} The first term $\chi: \prod_i \R^{n_i} \to [0,\infty)$ is a marginal cost function meant to penalize for travel-time spent away from some specified collection $\{\B x_{i,f}\}$ of desired final locations. It should, in some manner, be an approximation of the indicator function which gives value $0$ if $\B x_{i} = \B x_{i,f}$ for each $i$, and gives value $1$ if there is any $i$ for which $x_i \neq x_{i,f}$. In this way, $\chi$ encourages the agents to reach their final locations as quickly as possible to avoid accruing cost. The other term $\rho: \prod_i \R^{n_i} \to [0,\infty)$ is a marginal cost function meant to penalize the agents for breaking formation. This will be highly problem-dependent but should be approximately zero when the agents are approximately in the correct formation, and grow large when the agents are farther from being in the correct formation. The constants $w_1,w_2 > 0$ are weights that the controller can choose so as to prioritize time optimality or pattern formation.  

Fixing a point $\{x_i\}_{i = 1}^I$ in the holistic state-space. and $t \in [0,T]$, we define the value function \begin{equation} \label{eq:valueFunc} u(\{x_i\},t) = \inf_{\{\B a_i(\cdot)\}} \mathcal C_{\{x_i\},t}[\{\B x_i\}, \{\B a_i\}],  \end{equation} where $\mathcal C_{\{x_i\},t}$ is the cost functional \eqref{eq:CostFunc} restricted to the time interval $(t,T]$ and to trajectories $\{\B x_i\}$ such that $\B x_i(t) = x_i$ for each $i = 1,\ldots, I$.  In words: $u(\{x_i\},t)$ is the minimal remaining cost that the controller can achieve if the system is at state $\{x_i\}$ at time $t$. Our optimization problem is then equivalent to resolving $u(\{x_{i,0}\},0)$ along with the control values that achieved this minimal cost. Embedding this single problem within a whole class of similar problems allows one to use the dynamic programming principle \cite{Bellman} to conclude that, under mild conditions on the data (here $f_i, \chi, \rho$), the value function is formally the solution of the backward-in-time Hamilton-Jacobi-Bellman equation \begin{equation} \label{eq:HJB0}
\begin{split}
&u_t + \inf_{\{a_i\}}\left\{\sum^I_{i=1} \langle f_i, \nabla_i u\rangle\right\} + w_1\chi + w_2 \rho= 0,\\
&u(x,T) = 0.
\end{split} 
\end{equation} Here $\nabla_i$ denotes the gradient with respect to the variables describing the $i^{\text{th}}$ state-space, and we have suppressed the dependence of $f_i,\chi,\rho$ on the state variables for brevity. The general derivation as well as a discussion of so-called viscosity solutions for Hamilton-Jacobi equations can be found in several sources \cite{Bardi, Fleming, Liberzon}. Because forward-in-time PDE are more comfortable, it is relatively common to make the time reversing substitution $t \mapsto T-t$ to arrive at the equivalent forward in time equation \begin{equation} \label{eq:HJB} \begin{split} &u_t + \mathcal H\Big(\{x_i\}, \{\nabla_i u\}, t\Big) = 0,  \\ &u(x,0) = 0, \end{split}\end{equation} where the Hamiltonian $\mathcal H: \prod_i \R^{n_i} \times \prod_i \R^{n_i} \times [0,T] \to \R$ is given by 
\begin{equation} 
\label{eq:Hamiltonian}
\mathcal H\Big(\{x_i\},\{p_i\},t\Big) = \sup_{\{a_i\}} \left\{\sum^I_{i=1} \innerprod{-f_i}{p_i}\right\} -  w_1 \chi- w_2 \rho,
\end{equation} where we have used $p_i$ as a proxy for $\nabla_i u$. Assuming the solution of \eqref{eq:HJB} is known, the optimal feedback controls are given as the arguments that achieve the maximum in \eqref{eq:Hamiltonian}, and one can generate optimal trajectories by inserting the feedback controls into \eqref{eq:dynamics} and integrating the equations. 

We note that \eqref{eq:HJB} couples the behavior of the agents through the penalty terms $\chi$ and $\rho$, but decouples the agents' motion and control in the Hamiltonian, so that we can define individual Hamiltonians \begin{equation} \label{eq:HamiltonianIndiv}
H_i(x_i,p_i,t) = \sup_{a_i} \innerprod{-f_i}{p_i}
\end{equation} and re-write the PDE in \eqref{eq:HJB} as \begin{equation}
\label{eq:HJBsplit}
u_t + \sum^I_{i=1} H_i(x_i,\nabla_i u,t_i) = w_1 \chi\Big(\{x_i\}\Big) + w_2 \rho\Big(\{x_i\}\Big).
\end{equation}  We use \eqref{eq:HJB} when describing our numerical methods generally, but exploit the partial decoupling of the agents afforded by \eqref{eq:HJBsplit} in our implementation. 

In our demonstrations below, we we use the simple cases of isotropic motion as in \cite{Gee,ParkinsonPolage} and the Reeds-Shepp car \cite{ReedsShepp}. For isotropic motion, the dynamics are \begin{equation} \label{eq:isotropic}
\dot {\B x} = v(\B x,t) \B a, \,\,\,\,\,\,\, \abs{\B a(\cdot)} \le 1.
\end{equation} This could model a human walker or omnidirectional robot who can travel in any direction with local speed function $v(x,t)>0$. For isotropic dynamics, the individual Hamiltonian is \begin{equation}
\label{eq:isotropicH} H_i(x_i,p_i,t) = \sup_{\abs{a_i}\le 1} \innerprod{-v(x_i,t)a_i}{p_i} = v(x_i,t)\abs{p_i}. 
\end{equation} By contrast, the Reeds-Shepp car is a simplified model of vehicular motion wherein all motion must be tangential to the wheels and there is a minimum turning radius (or equivalently, a maximum angular velocity). This leads to curvature constrained paths modeled by the dynamics \begin{equation} \label{eq:RScar}
\dot {\B x}^{(1)} = v V\cos {\B x}^{(3)}, \,\,\,\,\, \dot {\B x}^{(2)} = vV\sin {\B x}^{(3)}, \,\,\,\,\, \dot {\B x}^{(3)} = W\omega.
\end{equation} Here the vehicle's configuration is given by $\B x = (\B x^{(1)}, \B x^{(2)}, \B x^{(3)})$ where the first two entries are the Cartesian coordinates for the center of mass, and the third entry is the angle describing orientation of the wheels \cite{TakeiTsai2,ParkinsonBoyle}. The normalized controls are $v(\cdot),\omega(\cdot) \in [-1,1]$, representing tangential and angular velocity respectively, and $V, W$ are constants denoting the maximum tangential and angular velocity respectively. Following through the derivation, the individual Hamiltonian for a Reeds-Shepp car is \begin{equation} \label{eq:RSH}
H_i(x_i,p_i,t) = V\vert p_i^{(1)}\cos x_i^{(3)} + p_i^{(2)}\sin x_i^{(3)}\vert  + W \vert p_i^{(3)}\vert
\end{equation} as derived in \cite{ParkinsonBoyle,TakeiTsai2}. Crucially, for both of these cases, the Hamiltonian can be explicitly resolved, which allows for seamless application of the numerical methods in the ensuing section. Investigation into the applicability of the methods when the Hamiltonian cannot be explicitly resolved, as in \cite{Parkinson}, is ongoing.

\section{NUMERICAL METHODS}

In this section, we describe the numerical methods used to approximately solve \eqref{eq:HJB}, as well as several practical implementation notes. The strategy is to approximately express the solution $u(\{x_i\},t)$ of \eqref{eq:HJB} as the solution of saddle point problem with regard to discrete state and co-state trajectories. We hint at these concepts here, while directing the reader to \cite{Lin,ParkinsonBoyle,Huynh} for a more full exposition.

Given a point $(\{x_i\},t)$ at which to resolve the solution of \eqref{eq:HJB}, we let $t_0<t_1<\ldots<t_J$ be a uniform discretization of $[0,t]$ with grid spacing $\delta = t/J,$ and for each agent $\B x_i$, let $x_{i,j}$ be an approximation to the value $\B x_i(t_j)$. Likewise, we introduce a discrete costate trajectory $p_{i,j}$ for each agent. These are the values of $\nabla_i u(t_j)$ along the trajectory, and can be used to compute the optimal control along the trajectory as the minimizer in \eqref{eq:Hamiltonian}. Mathematically, they can be seen to arise in \eqref{eq:valueFunc} if one discretizes the cost functional and enforces the dynamics \eqref{eq:dynamics} using Lagrange multipliers \cite{Lin,ParkinsonBoyle}. We use $\{x_{i,j}\} \subset \prod_i \R^{n_i}$ to denote the approximate position of all $i$ agents at time step $j$. The time-reversing substitution for the HJB equation needs to be carried through here, so these trajectories will be time-reversed optimal trajectories.

 Using this notation, the authors of \cite{Lin,ParkinsonBoyle,Huynh} conjecture that the solution of \eqref{eq:HJB} can be approximated by \begin{equation}
\label{eq:HopfLaxType}
\begin{split} u(\{x_i\},t) \approx 
\inf_{\{x^*_{i,j}\}} \sup_{\{p^*_{i,j}\}}\bigg\{ &\sum_{j} \sum_i \innerprod{p^*_{i,j}}{x^*_{i,j} - x^*_{i,j-1}} \\ -  &\delta \sum_j\mathcal H\Big(\{x^*_{i,j}\},\{p^*_{i,j}\},t_j\Big) \bigg\}. \end{split} 
\end{equation} The values at the final point $\{x_{i,J}\}$ should be set to  $\{x_i\}$, the point at which one desires to solve \eqref{eq:HJB}. Solving this saddle-point problems provides the value $u(\{x_i\},t)$ as well as the optimal state and costate trajectories. 

We point out the formal similarity between \eqref{eq:HopfLaxType} and \eqref{eq:ClassicHL}. If we stack all agents into a single vector $x$, we can write the sums over $i$ in both \eqref{eq:Hamiltonian} and \eqref{eq:HopfLaxType} as a inner product in a higher dimensional space. Then accounting for the zero initial condition in \eqref{eq:HJB}, \eqref{eq:HopfLaxType} and \eqref{eq:ClassicHL} only differ in that the time interval $[0,t]$ has been discretized and we consider points along the optimal trajectory in \eqref{eq:HopfLaxType}. If the optimal costate trajectory happens to be a straight line (as is the case when the Hamiltonian is time- and space-independent like in \eqref{eq:ClassicHL}), these sums telescope and one can recover \eqref{eq:ClassicHL} from \eqref{eq:HopfLaxType}. Accordingly, one could see \eqref{eq:HopfLaxType} as a generalized, discrete Hopf-Lax formula. It bears repeating that at this point, it is merely conjectural that the saddle-point problem in \eqref{eq:HopfLaxType} does indeed approximate the value function, though solid numerical evidence is given in \cite{ParkinsonPolage,Lin,ParkinsonBoyle,Huynh} for a variety of Hamiltonians, and our results corroborate this evidence. 

We solve the saddle point problem \eqref{eq:HopfLaxType} using a version of the Chambolle-Pock Primal Dual Hybrid Gradient algorithm \cite{PDHG} adapted to our particular problem. The full algorithm is detailed in Algorithm \ref{alg:1}. It is an iterative algorithm which alternately performs the optimization over state and costate variables. A key realization is that the optimization in \eqref{eq:HopfLaxType} is at least partially decoupled in both $i$ and $j$. The decoupling allows us to deal solely with the values for individual agents at individual time steps, using the individual Hamiltonians, rather than considering the whole trajectory for the entire group of agents all at once. This is the reason that $\mathcal H$ does not appear in Algorithm \ref{alg:1}. Any optimization that involves $\mathcal H$ is broken down into parts that only involve $H_i$. In the algorithm, we use $\chi_i$ and $\rho_i$ to denote the functions $\chi$ and $\rho$ treated as functions of $x_i$ only, while all other arguments are fixed at their current values. Note that these functions depend only on state variables, and are thus omitted from the costate optimization step. 

We include some practical implementation notes for Algorithm \ref{alg:1}. A first note is that there has been no discretization of spatial domains, only of time, so Algorithm \ref{alg:1} can resolve the solution of \eqref{eq:HJB} at single points (and thus compute optimal trajectories) in a manner that scales well to high dimensions, which is crucial for multi-agent problems. However, in the ensuing examples, several thousand iterations are required for the algorithm to converge, so it is of paramount importance that each iteration be computed as efficiently as possible. To this end, it behooves one to evaluate the minimization problems in Algorithm \ref{alg:1} explicitly whenever possible. For example, the case of isotropic motion or Reeds-Shepp cars, where the individual Hamiltonians are given by \eqref{eq:isotropicH} and \eqref{eq:RSH} respectively, it is possible to resolve the minimization for the costate variables explcitly. For an agent exhibiting isotropic motion \eqref{eq:isotropic}, the update is \begin{equation} \label{eq:pupdateIsotropic}
p_{i,j} \leftarrow \max\left\{0,1-\tfrac{\sigma \delta v(x_{i,j})}{\abs{\beta}}\right\}\beta,
\end{equation} and for Reeds-Shepp dynamics \eqref{eq:RScar}, the update is \begin{equation}
\label{eq:pupdateRS}
\begin{split}
p^{(1:2)}_{i,j} &\leftarrow \beta^{(1:2)} -\min\left\{ 1, \tfrac{V\sigma \delta}{\abs{\gamma^T\beta^{(1:2)}}} \right\}(\gamma^T\beta^{(1:2)})\gamma, \\
p^{(3)}_{i,j} &\leftarrow \max \left\{0, 1-\tfrac{W\sigma \delta }{\beta^{(3)}}\right\}\beta^{(3)},
\end{split}
\end{equation} where in both cases $\beta$ is as in Algorithm \ref{alg:1}, and in \eqref{eq:pupdateRS},  $\gamma = (\cos x_{i,j}^{(3)}, \sin x_{i,j}^{(3)})$, and the superscript $(1:2)$ denotes the first two coordinates of the vector. 

\begin{algorithm}[t!]
\caption{Algorithm for Solving \eqref{eq:HopfLaxType}}
Input the point $(\{x_i\},t)$ at which to resolve the HJB equation and hyperparameters $\sigma, \tau > 0$ with $\sigma\tau \le 0.25$. 

 Set $\{x_{i,J}\}= \{x_i\} , \{p_{i,0}\} = \{0\}$, and initialize $\{x_{i,j}\}$, $\{p_{i,j}\}$ randomly for all other values. Set $\{z_{i,j}\} = \{x_{i,j}\}$ for all $i,j$. 

\begin{algorithmic}[t!]
\Repeat
 \For {$j = 1$ to $J$}
    	\For {$i = 1$ to $I$}
    		\State $\beta \leftarrow p_{i,j}^k + \sigma(z_{i,j}^k - z_{i,j-1}^k)$
    		\State $p_{i,j} \leftarrow \argmin_{\tilde p \in \R^{n_i}} \{ \delta H_i(t_j,x_{i,j},\tilde p) + \frac{1}{2\sigma} \lvert \tilde p- \beta \rvert^2 \}$
	\EndFor
    \EndFor
    
    \State $x^{\text{old}}_{i,j} \leftarrow x_{i,j}$ \,\,\,\,\,\,\,\,\,\,\,\,\,\,\,\,\,\,\, (for each $i,j$)
    \State $x_{i,0} \leftarrow x_{i,0} + \tau p_{i,1}$\,\,\,\,\,  (for each $i$)
    \For{$j = 1$ to $J - 1$}
    \For {$i = 1$ to $I$ }
    \State $\xi \leftarrow x_{i,j} - \tau(p_{i,j} - p_{i,j+1})$
    \State $x_{i,j} \leftarrow \argmin_{\tilde x \in \R^{n_i}}\{ - \delta H_i(t_j,\tilde x,p_{i,j}) - w_1 \delta \chi_i(\tilde x)$ 
    \State $\tfrac{}{}$ \hspace{2.5cm} $- w_2 \delta \rho_i(\tilde x)+ \frac{1}{2\tau} \lvert \tilde x -\xi \rvert^2 \}$
    \EndFor
    \EndFor
    
    \State $z_{i,j} \leftarrow 2x_{i,j} -x^{\text{old}}_{i,j}$  \,\,\,\,\, (for each $i,j$)
    
\Until{convergence}

\State $u =\sum_{i,j} \langle p_{i,j}, x_{i,j} - x_{i,j-1} \rangle - \delta \sum_j \mathcal H(t_j,\{x_{i,j}\},\{p_{i,j}\})$
\State \textbf{return } $u$, $\{x_{i,j}\}$
\end{algorithmic}
\label{alg:1}
\end{algorithm}

Except in extremely simple cases, it will \emph{not} be possible to explicitly solve the optimization problems to update the state variables in Algorithm \ref{alg:1}. In this case, one must approximate the optimizers. It is observed in \cite{Lin} that this approximation can actually be quite crude and the algorithm will still work. Following the suggestion of \cite{Lin}, we use a single step of gradient descent with some predefined step size $\eta$. As in \cite{ParkinsonPolage}, we begin with a large gradient descent rate and, after 5000 iterations, we decrease the gradient descent step size every 1000 iterations to resolve the path on a finer scale. 

As noted above, the choice of $\rho$, the marginal penalty function for breaking formation, is highly problem dependent, so we specify $\rho$ individually for the examples below. For $\chi$, the marginal penalty function for spending time away from the desired final points, we use \begin{equation} \label{eq:chidef}\chi(\{x_i\}) = \sum_i \Big(1-e^{-\abs{x_i-x_{i,f}}^2}\Big). \end{equation} Here we are treating $1-e^{-\abs{x_i - x_{i,f}}^2}$ as a smooth approximation to the indicator function which is $1$ when $x_i \neq x_{i,f}$ and $0$ when $x_i = x_{i,f}$. We could normalize $\chi$ by dividing by the total number of agents, so that it takes values in $[0,1]$, but because it is always appended with the weight $w_1$, the normalization is somewhat arbitrary. An advantage of using this particular form for $\chi$ is that it also decouples the agents, so that the only coupling comes in the formation function $\rho$. 

One last concern is how to add impassible obstacles computationally.We include obstacles in a manner similar to that proposed in \cite{TakeiTsai2,ParkinsonPolage,ParkinsonBoyle}. We set the velocity for all agents to zero inside obstacles. In this manner, if an agent enters an obstacle, they become stuck and cannot make it to their final point, and will incur a large cost. In practice, if $\Omega(t)$ is the subset of the physical domain occupied by obstacles at time $t$, we define the obstacle function \begin{equation} \label{eq:obstacles} 
O(x,t) = \tfrac 1 2 (1+\tanh(100d(x,t)^2)
\end{equation} where $d(x,t)$ is the signed distance from $x$ to the set of obstacles at time $t$ (which has negative values inside the obstacles). This function is approximately 0 inside obstacles and approximately 1 outside the obstacles. This introduces the issue of being able to compute $d(x,t)$ quickly. For this, we use the strategy of \cite{ParkinsonBoyle} wherein all obstacles are required to be circles, and the distance function is computed from simple geometric considerations. More complicated obstacles are approximated by collections of disjoint circles if desired. In implementation, any velocity term (for example $v(x,t)$ in \eqref{eq:isotropicH}, and $V,W$ in \eqref{eq:RSH}) is simply multiplied by $O(x,t)$. We use the smooth approximation of $O(x,t)$ because these terms will participate in the gradient descent which is used to resolve the updated state variable at each iteration in Algorithm \ref{alg:1}.

\section{RESULTS AND DISCUSSION}

To conclude, we present some examples demonstrating the efficacy of our algorithm. In each case, the agents navigate in two spatial dimensions, and the hyperparameters for Algorithm \ref{alg:1} are $\sigma = 1$, $\tau = 0.25$. We repeat the iteration until successive iterates do not differ in supremum norm by more than $5\times 10^{-4}$. For scale reference, each of these examples takes place roughly within the domain $[-2,2]^2$. The examples included are sythentic and are meant to demonstrate that multi-agent path planning with nonlinear dynamics can be efficiently computed using a Hamilton-Jacobi-Bellman formulation and the generalized Hopf-Lax formula \eqref{eq:HopfLaxType}. All simulations were run in MATLAB on the first author's desktop computer with an Intel(R) Core(TM) i7-14700 processor running at 2100Mhz with 32GB of RAM. The code that generated these examples (as well as animated versions of these examples) can be found at the link below.\footnote{\href{https://github.com/chparkinson/HL_multi_agent_pattern/}{\texttt{github.com/chparkinson/HL\_multi\_agent\_pattern/}}}

\begin{figure*}[!]
\centering
\fbox{\includegraphics[width=0.16\textwidth,trim = 80 25 65 20, clip]{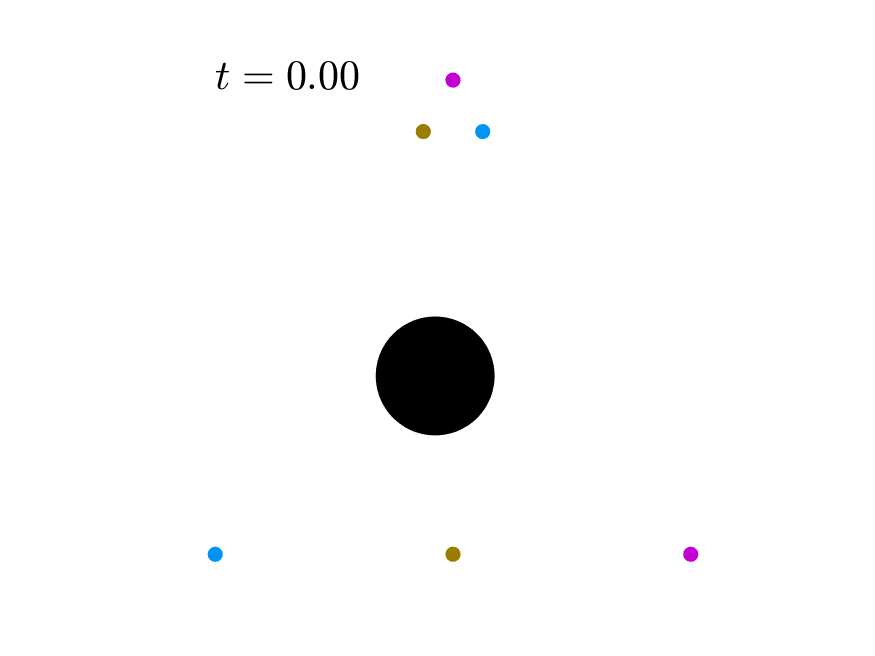}} \, 
\fbox{\includegraphics[width=0.16\textwidth,trim = 80 25 65 20, clip]{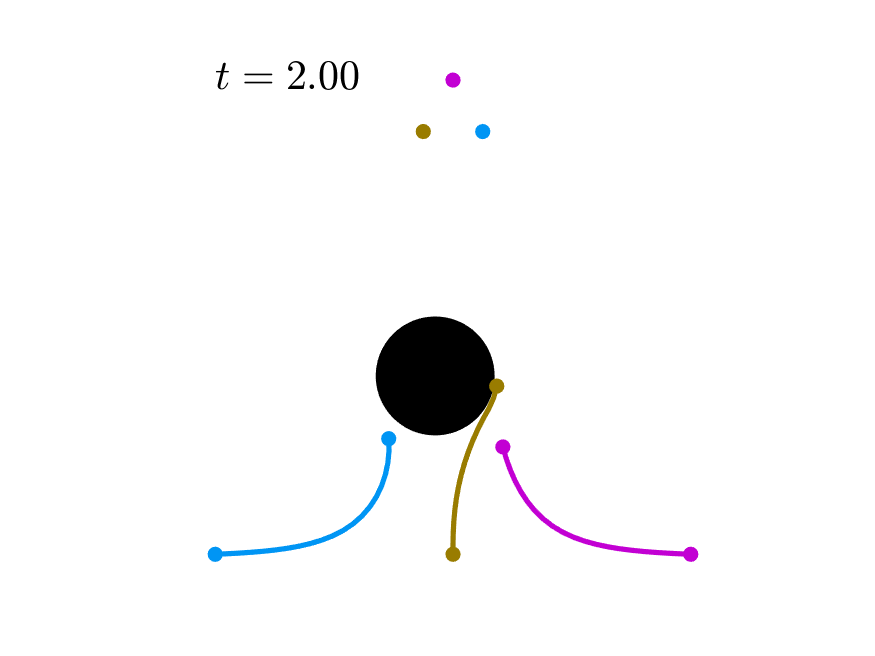}} \, 
\fbox{\includegraphics[width=0.16\textwidth,trim = 80 25 65 20, clip]{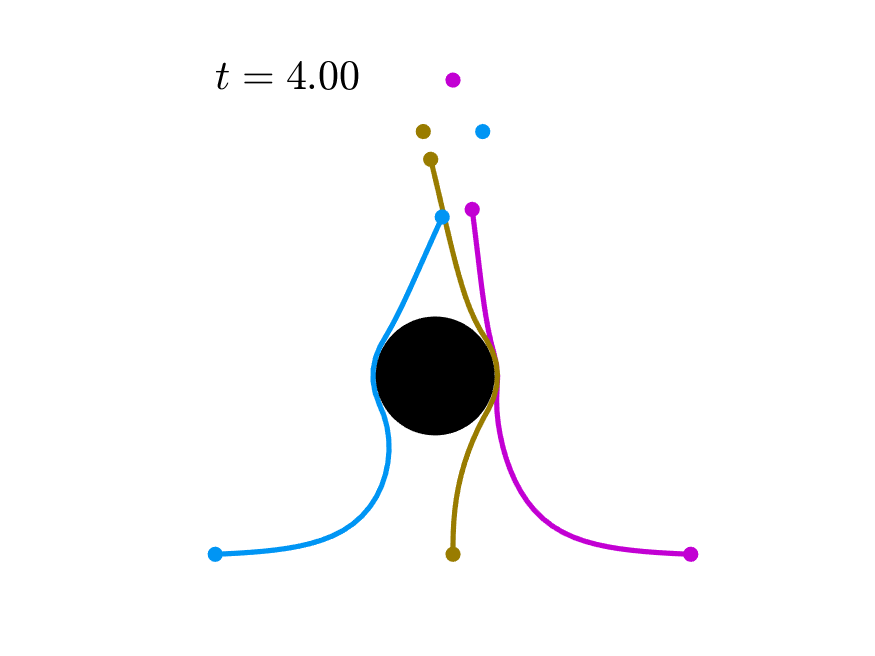}}\, 
\fbox{\includegraphics[width=0.16\textwidth,trim = 80 25 65 20, clip]{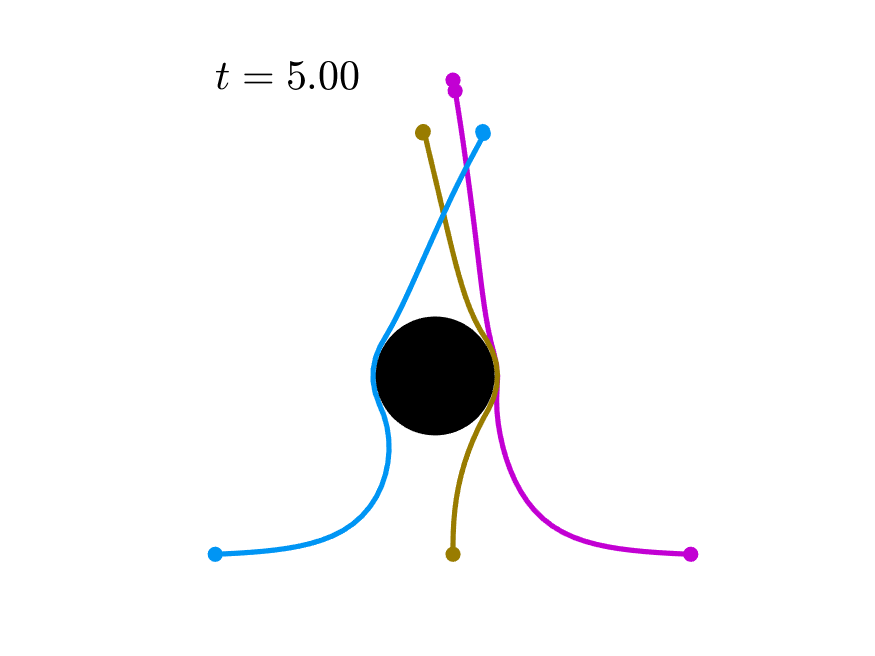}} \, 
\fbox{\includegraphics[width=0.16\textwidth,trim = 80 25 65 20, clip]{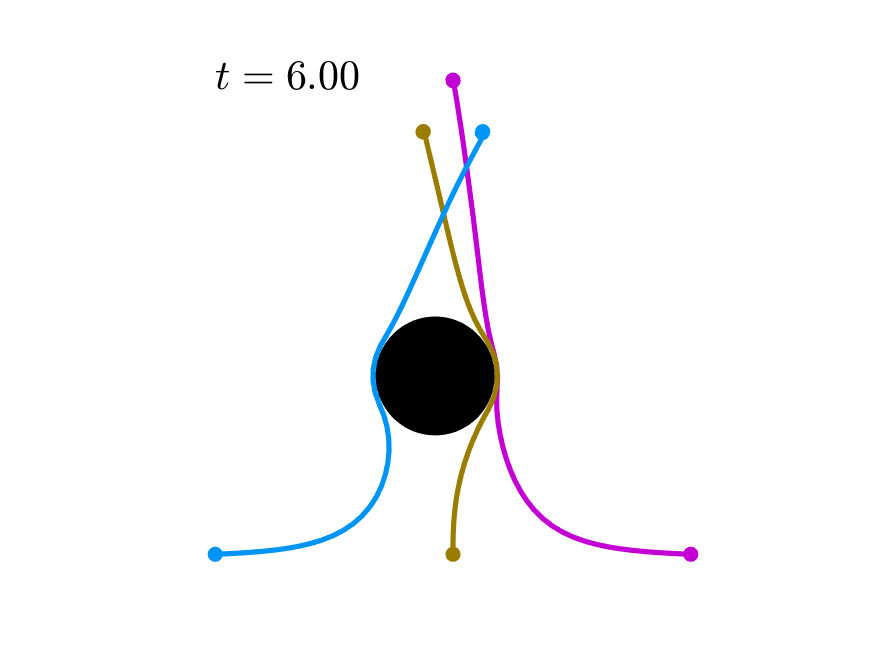}} \\ 
\vspace{3mm}
\fbox{\includegraphics[width=0.16\textwidth,trim = 80 25 65 20, clip]{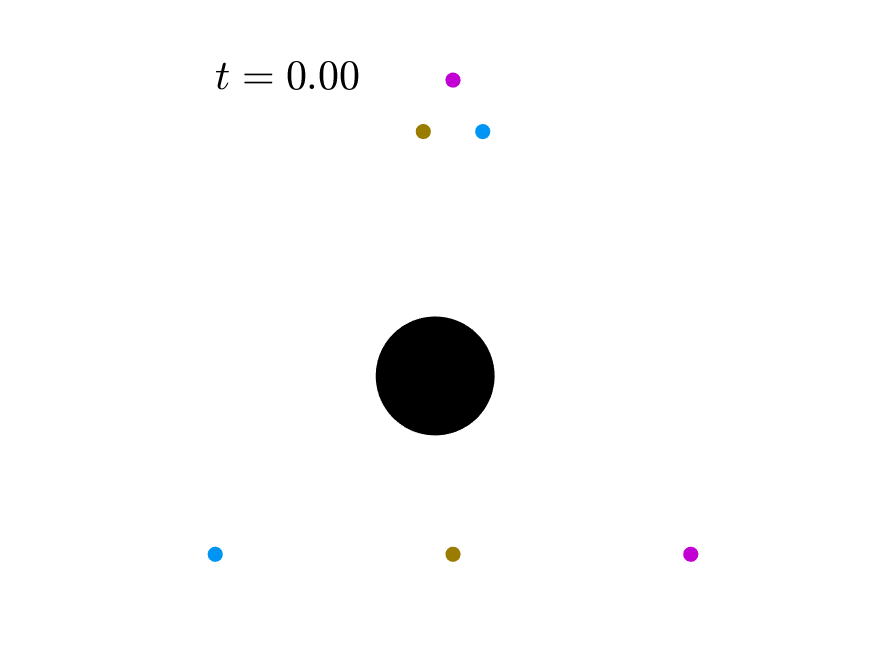}} \, 
\fbox{\includegraphics[width=0.16\textwidth,trim = 80 25 65 20, clip]{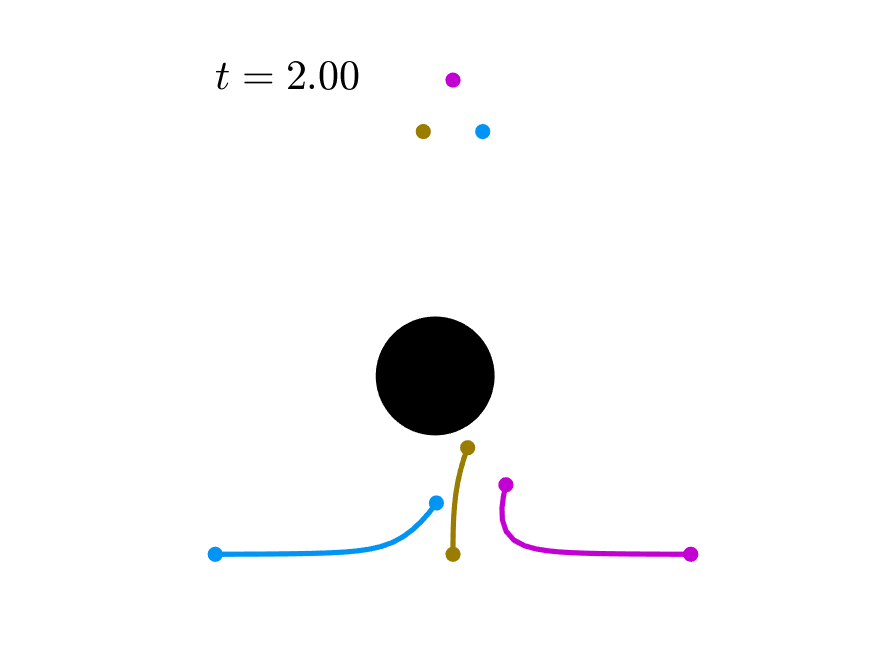}} \, 
\fbox{\includegraphics[width=0.16\textwidth,trim = 80 25 65 20, clip]{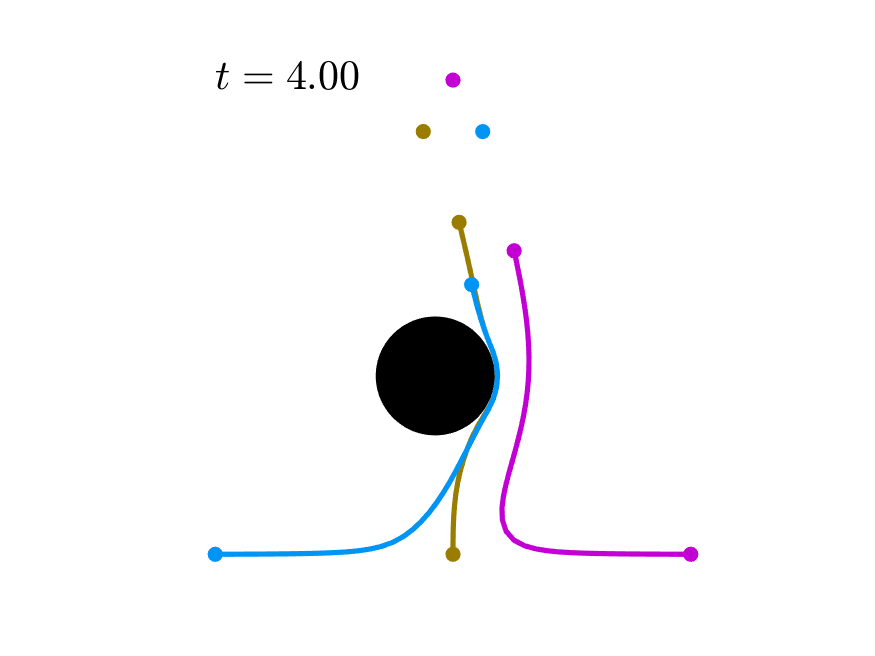}}\, 
\fbox{\includegraphics[width=0.16\textwidth,trim = 80 25 65 20, clip]{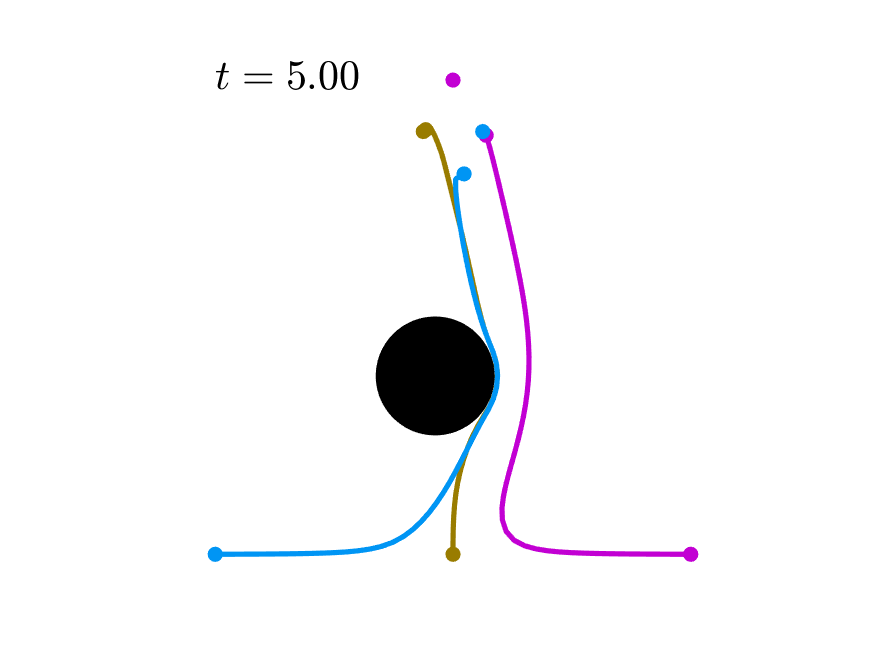}} \, 
\fbox{\includegraphics[width=0.16\textwidth,trim = 80 25 65 20, clip]{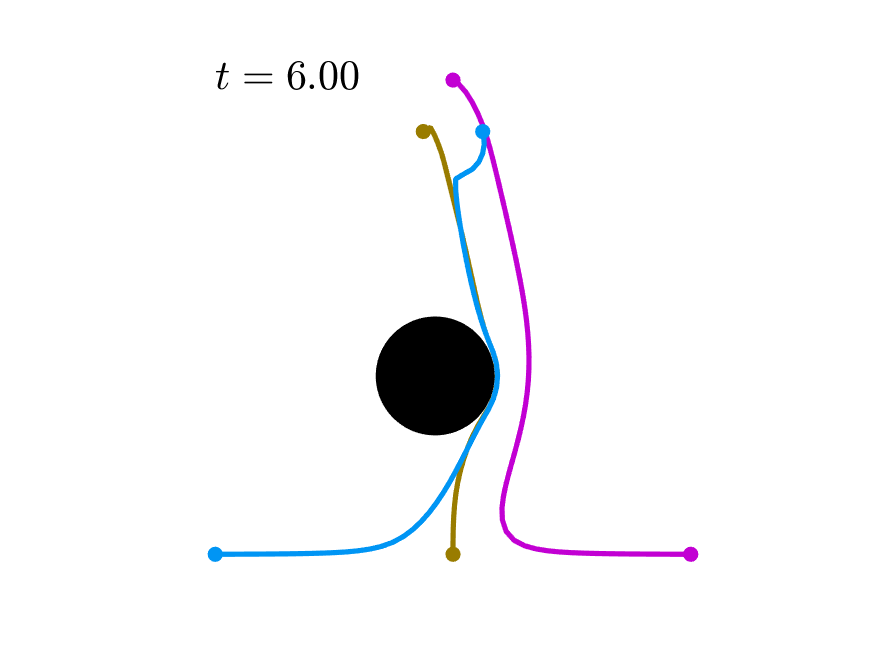}} 
\caption{Three agents navigate around an obstacle while trying to maintain formation in an equilateral triangle. Top: the agents are heavily incentivized to minimize travel time to their final destinations. Bottom: the agents are heavily incentivized to maintain formation. }
\label{fig:1}
\end{figure*}

In our first example in figure \ref{fig:1}, three agents (blue, yellow, purple dots) exhibiting isotropic motion with constant velocity $1$ travel from the bottom the frame to the top of the frame while attempting to form an equilateral triangle of side length 1/2. With 3 agents each each occupying a 2-dimensional state space, the total state space in this example is $6$ dimensional. In this case the marginal penalty function is $\rho(\{x_i\}) = \sum_{i\neq j}(\vert x_i - x_j\vert^2-0.25)^2$. The desired final points are at the top of each frame and do indeed form such a triangle. There is one impassible obstacle (black circle). In this case, we include simulations with two different values of the weights $w_1,w_2$. In the first simulation (top frames), we set $(w_1,w_2) = (1, \tfrac 1 2)$ indicating that it is more important to minimize travel time than it is to keep formation. Accordingly, the agents that begin in the corners do indeed bend inward a bit to get closer to forming a triangle, but ultimately, the agent on the left travels around the obstacle on the opposite side as the other two, and the agents reach their final points at roughly $t = 5$ seconds. By contrast in the second simulation (bottom frames), the weights are $(w_1,w_2) = (\tfrac 1 2, 4)$, indicating that it is very much important to keep formation. Accordingly, the middle agent does not start moving until the corners have drawn in, and then all agents travel the same direction around the obstacle in a triangular formation. Their dedication to keeping the formation is demonstrated in the final two frames where we see that, as they approach their final configuration, their orientation is off, so the first agent finds its final point, and the agents rotate while keeping formation to place the others at their final points. While there is randomness due to the random initialization, Algorithm \ref{alg:1} found optimal trajectories for this example in roughly 2 seconds of CPU time. 

\begin{figure*}[!]
\centering
\fbox{\includegraphics[width=0.16\textwidth,trim = 80 25 65 20, clip]{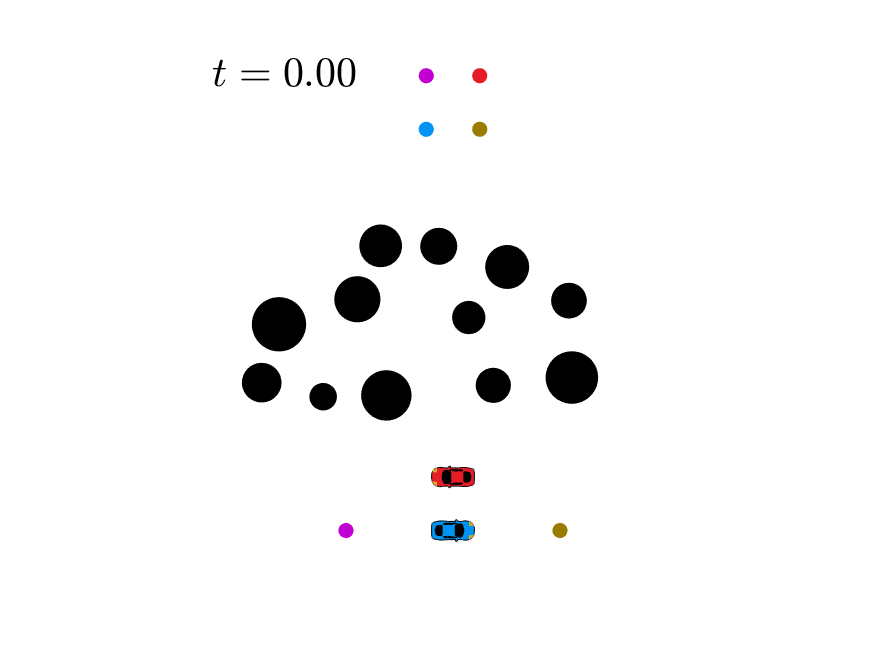}} \, 
\fbox{\includegraphics[width=0.16\textwidth,trim = 80 25 65 20, clip]{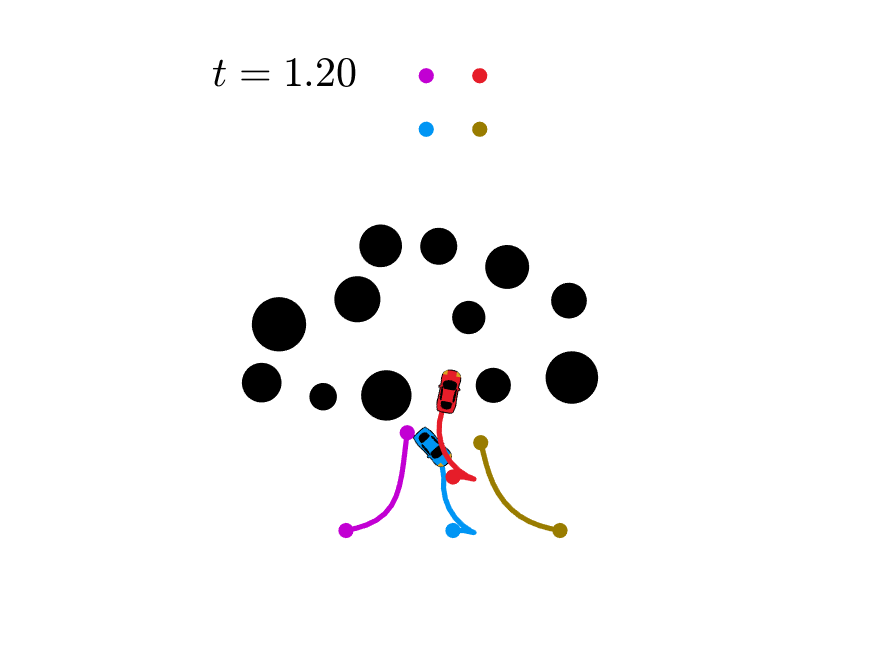}} \, 
\fbox{\includegraphics[width=0.16\textwidth,trim = 80 25 65 20, clip]{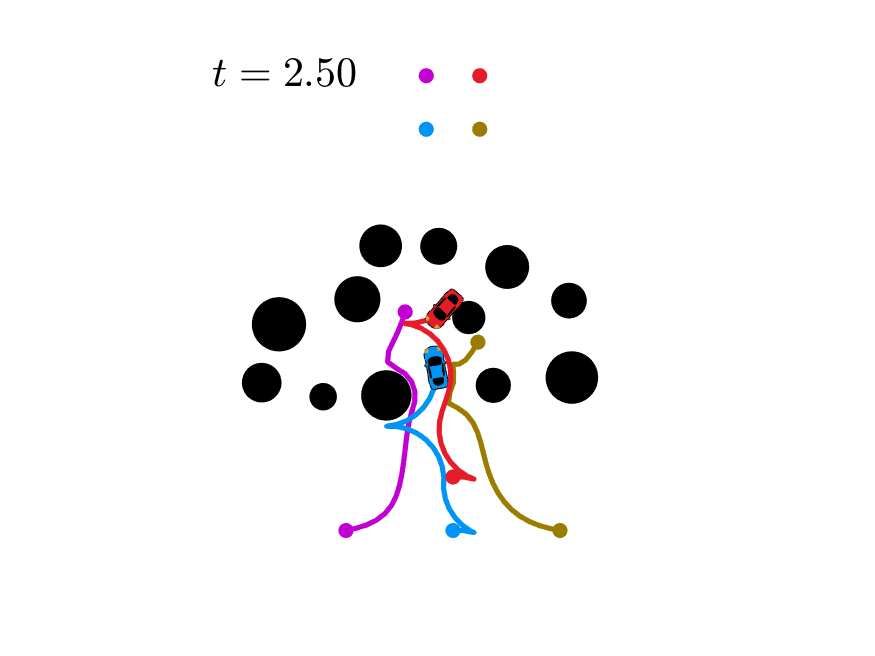}}\, 
\fbox{\includegraphics[width=0.16\textwidth,trim = 80 25 65 20, clip]{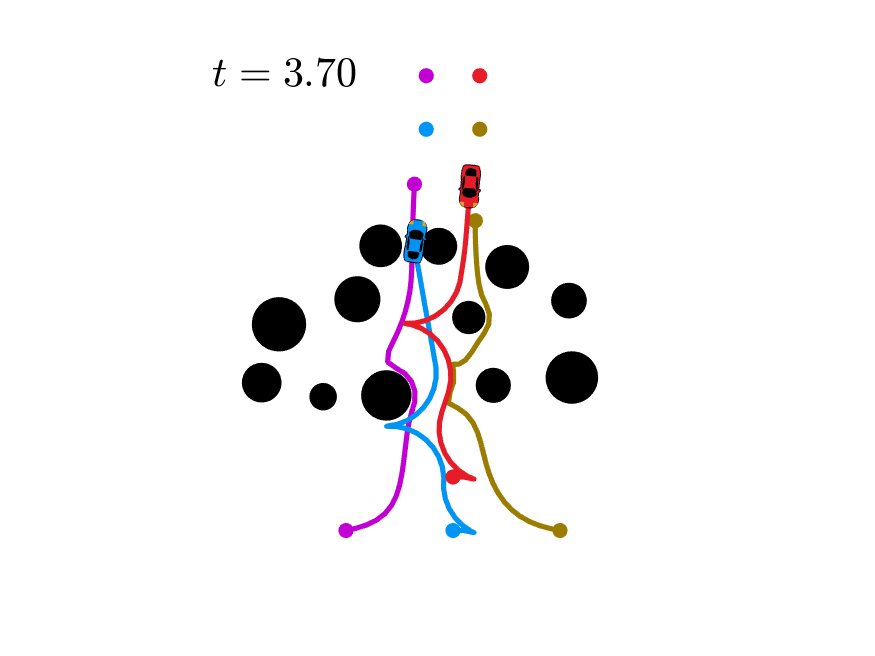}} \, 
\fbox{\includegraphics[width=0.16\textwidth,trim = 80 25 65 20, clip]{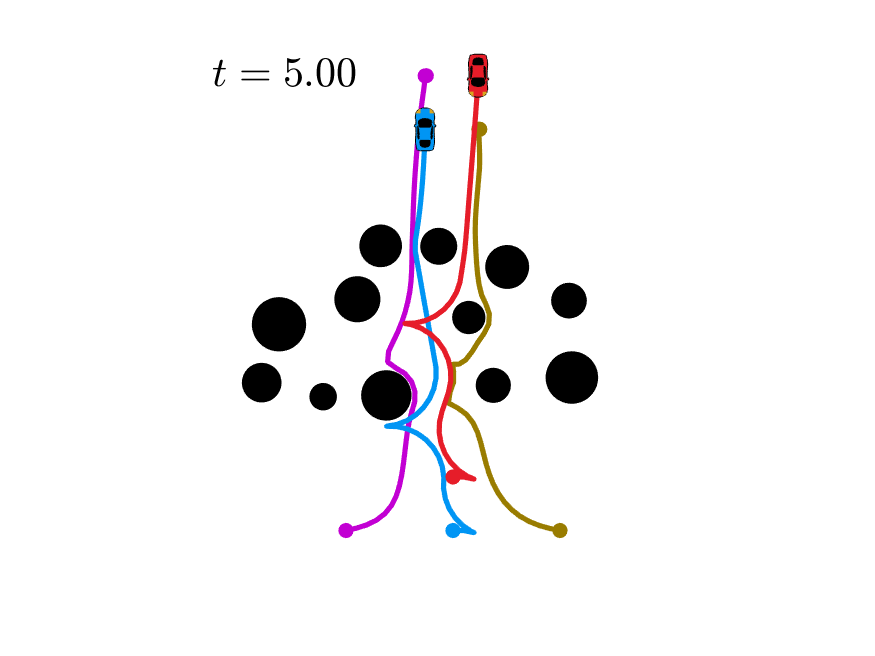}} \\ 
\caption{Four agents navigate around obstacles while trying to maintain formation in a square. The two agents exhibiting isotropic motion (dots) mimic the cuspidal paths that the cars take so as to keep in formation, even though their dynamics allow for smooth paths.}
\label{fig:2}
\end{figure*}

In a second example, we have four agents, two of which exhibit isotropic motion with velocity 1 (yellow and purple dots) and two of which have Reeds-Shepp dynamics (red and blue cars) with maximum tangential and angular velocities $V = 1$, $W = 2$. They navigate through a field of obstacles while maintaining formation in a square of side length 1/2,  where the cars comprise opposite corners. Results of this simulation are in figure \ref{fig:2} This demonstrates the effectiveness of our method in handling heterogeneous agents. Four points $A,B,C,D$ in the plane form a square of side length $d$ if and only if $\abs{A-C}^2 - 2d^2 = \abs{B-D}-2d^2 = 0, \abs{(A+C)-(B+D)}=0$, and $\innerprod{A-C}{B-D} = 0$. Accordingly our penalty function is the sum of the squares of each of these formulas, where $A,C$ are the locations of the isotropic agents and $B,D$ are the locations of the cars, and we take $(w_1,w_2) = (1,1)$.  One noteworthy feature here is that Reeds-Shepp cars are known to travel along curvature constrained paths with cusps at the spots where they reverse direction \cite{ReedsShepp}. Here the isotropic agents mimic the cusps so as to keep formation, even though their dynamics allow them to travel along smooth paths. The collective dynamics take place in a 10 dimensional state space, and this simulation required on the order of 10 seconds of CPU time to resolve. 

\begin{figure*}[!]
\centering
\fbox{\includegraphics[width=0.16\textwidth,trim = 90 35 75 25, clip]{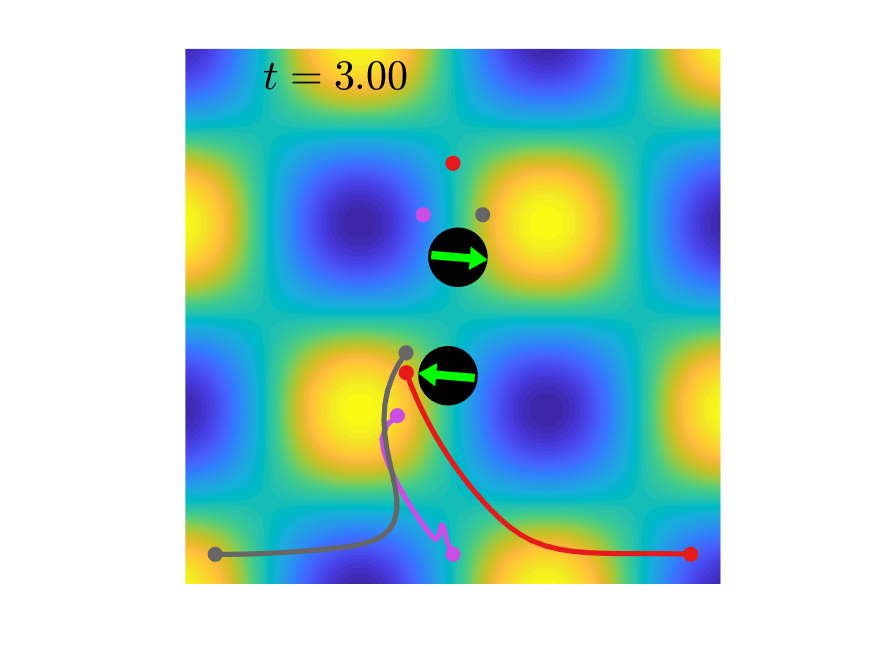}} \, 
\fbox{\includegraphics[width=0.16\textwidth,trim = 90 35 75 25, clip]{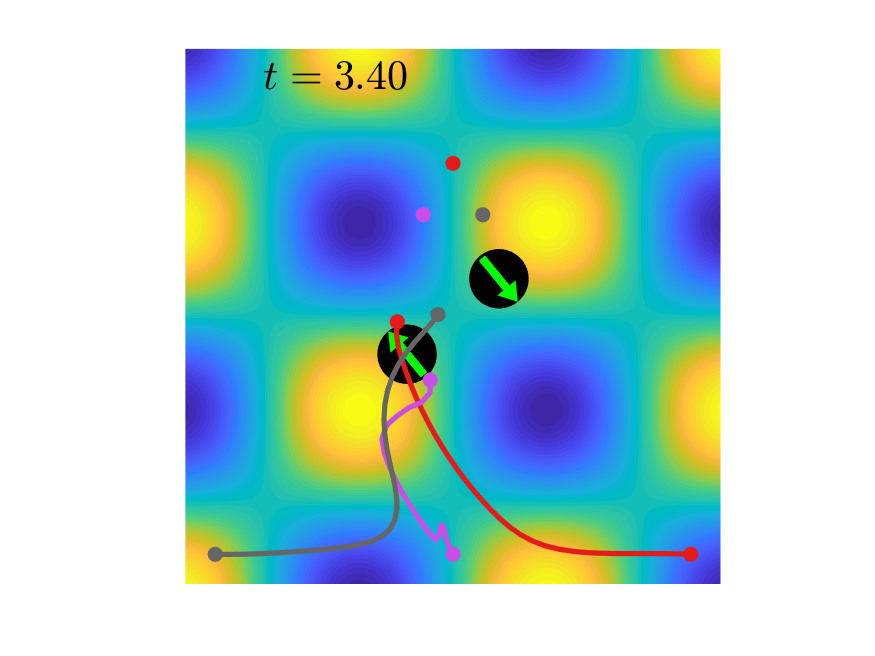}} \, 
\fbox{\includegraphics[width=0.16\textwidth,trim = 90 35 75 25, clip]{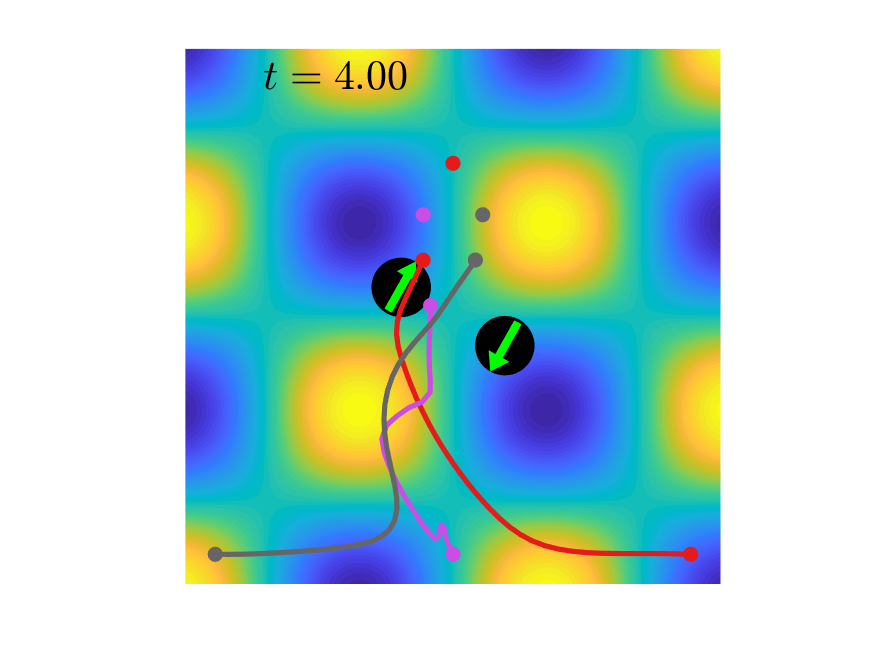}}\, 
\fbox{\includegraphics[width=0.16\textwidth,trim = 90 35 75 25, clip]{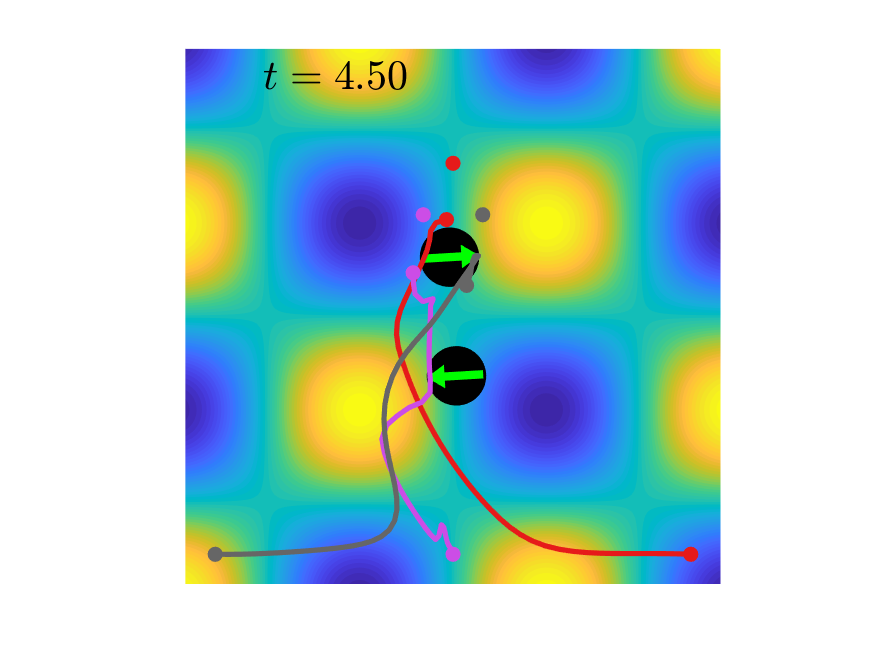}} \, 
\fbox{\includegraphics[width=0.16\textwidth,trim = 90 35 75 25, clip]{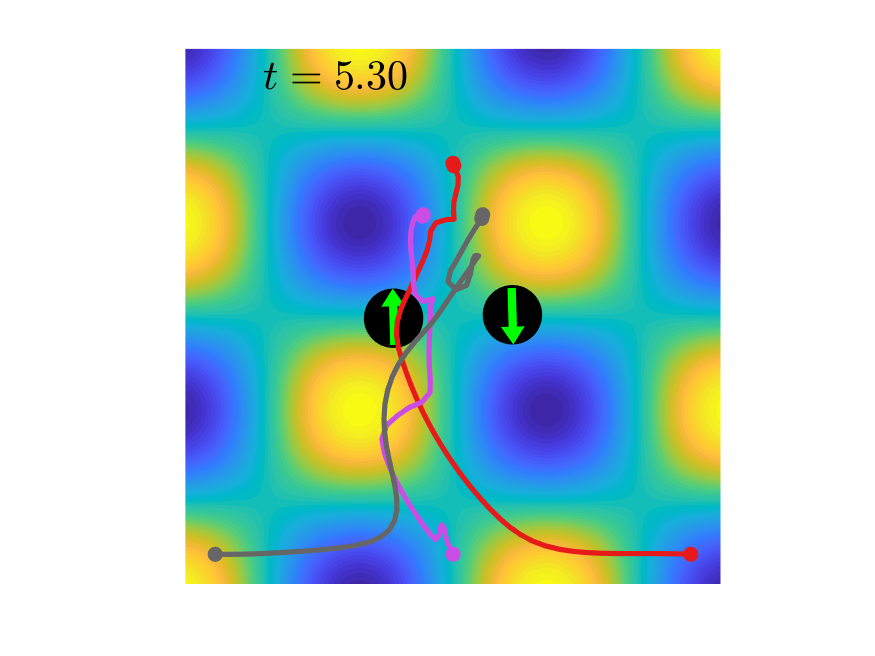}} \\ 
\caption{Three agents navigate around two moving obstacles while trying to maintain formation in an equilateral triangle.  They move with local speed $v(x,y) = 1+\frac 1 4 \sin(x)\sin(y)$ (faster in yellow regions; slower in blue). The obstacles are rotating in a circle with direction specified by the red arrow. We have chosen specifically timed snapshots to demonstrate interaction with the obstacles.}
\label{fig:3}
\end{figure*}

In a third and final example, we have the same setup as in figure \ref{fig:1} where three agents exhibit isotropic motion and try to maintain a triangle. However, in this case, there is an underlying speed function $v(x,y) = 1+\frac 1 4 \sin(x)\sin(y)$ so that the agents can move faster in some regions and slower in others. Additionally, there are two obstacles which rotate in a circle as time progresses, so that this is a time dependent problem, and the final points have been moved down closer to the obstacles. We take $(w_1,w_2) = (\frac 1 2,3)$ so that the agents are heavly incentivized to keep formation. The results are in figure \ref{fig:3}, where we display specifically timed snapshots to demonstrate the manner in which the agents interact with the moving obstacles. This demonstrates the utility of our algorithm on time dependent problems. These trajectories were typically resolved in roughly 4 seconds. 

\bibliographystyle{ieeetr}
\bibliography{Biblio}

\end{document}